\documentclass{gtart}
\usepackage{amscd,amssymb,amsxtra}

\input gtoutput
\volumenumber{2}\papernumber{8}\volumeyear{1998}
\pagenumbers{145}{174}\published{17 August 1998}
\proposed{Haynes Miller}
\seconded{Ralph Cohen, Gunnar Carlsson}
\received{5 September 1997}\revised{27 March 1998}
\accepted{17 August 1998}

\numberwithin{equation}{section}
\newtheorem{Theorem}[equation]{Theorem}

\newtheorem{Lemma}[equation]{Lemma}
\newtheorem{Proposition}[equation]{Proposition}
\newtheorem{Corollary}[equation]{Corollary}
\newtheorem{ITheorem}{Theorem}
\theoremstyle{definition}  
\newtheorem{Definition}[equation]{Definition}

\renewcommand{\th}{\mbox{$^{\text{th}}$ }}

\newcommand{\Tate}{t}
\newcommand{\eqdef}{\overset{\text{def}}{=}}
\newcommand{\EJW}[1]{E (#1)}     
\newcommand{\ELT}[1]{E_{#1}}     
\newcommand{\EFP}[1]{\EJW{#1}}   
\newcommand{\ERT}{\EJW{n}[\w]}   

\newcommand{\bp}[1]{BP\langle #1  \rangle}
\newcommand{\jw}[1]{E (#1)}
\newcommand{\cjw}[1]{E^{\wedge } (#1)}
\newcommand{\ltjw}[1]{E_{#1}}
\newcommand{\gltjw}[1]{\ltjw{#1}}
\newcommand{\Esource}{T E}

\newcommand{\psb}[1]{[ \! [#1] \! ]}
\renewcommand{\k}{\mathbf{F}_{p}}
\newcommand{\F}{\mathbf{F}}

\newcommand{\C}{\mathbf{C}}
\newcommand{\SepClos}[1]{#1^{\textup{sep}}}
\newcommand{\khuge}{\SepClos{\k\lsb{\y}}}
\newcommand{\Iadic}[1]{^{\wedge}_{I_{#1}}}
\newcommand{\Iwadic}{^{\wedge}_{I_{n-1}}}
\newcommand{\tensor}[1]{\underset{#1}{\otimes}}

\newcommand{\fs}[1]{\underset{#1}{+}}
\newcommand{\fdiff}[1]{\underset{#1}{-}}
\newcommand{\lsb}[1]{( \! (#1) \! )}
\newcommand{\GJW}[1]{G (#1)}
\newcommand{\GLT}[1]{G_{#1}}
\newcommand{\GNice}{F}
\newcommand{\GNiceRes}{\GNice_{0}}
\newcommand{\E}{\GNiceRes}
\newcommand{\Gabbrev}{\overline{F}}
\newcommand{\Gtypical}{G}
\newcommand{\GH}{H}
\renewcommand{\H}{\GH}
\newcommand{\Honda}{\GH}

\newcommand{\rtvmo}{w}
\newcommand{\w}{\rtvmo}

\newcommand{\proottwo}{x}
\newcommand{\x}{\proottwo}

\newcommand{\prootzero}{y}
\newcommand{\y}{\prootzero}

\newcommand{\dps}{[p]}

\newcommand{\st}{\mbox{$\star$}}
\newcommand{\notdivide}{\nmid}

\newcommand{\iadiccomp}[2]{#1^{\wedge}_{I_{#2}}}
\newcommand{\pcomp}[1]{#1_{p}\!\!\!\hbox{\^{ }}}

\newcommand{\Zloc}{\ints_{(p)}}

\newcommand{\Zp}{\ints_{p}}
\newcommand{\zp}{\ints/ (p)}

\newcommand{\im}{\mbox{\rm im }}
\newcommand{\ints}{\mathbf{Z}}
\def\invlim#1#2{\mathchoice
        {\mathop{\displaystyle\mathop{\mbox{\rm lim}}_{\longleftarrow}}
        _{#1}{#2}}
        {\mathop{\smash{\displaystyle\mathop{\mbox{\rm lim}}
        _{\longleftarrow}}}
        _{#1}{#2}\vphantom{\displaystyle\mathop{\mbox{\rm lim}}
        _{\longleftarrow}}}
        {}{}}

\begin{document}

%
%
%
%
\title{Completions of $\mathbf{Z}/ (p)$--Tate cohomology of\\ periodic spectra}
\covertitle{Completions of ${\noexpand\bf Z}/(p)$--Tate cohomology of\\ periodic spectra}
\asciititle{Completions of Z/(p)-Tate cohomology of periodic spectra}
\shorttitle{Tate cohomology spectra}

\authors{Matthew Ando\\Jack Morava\\Hal Sadofsky}
\shortauthors{Ando, Morava and Sadofsky}

\address{Department of Mathematics,
Unversity of Virginia\\
Charlottesville, VA  22903\\\smallskip
\\Department of Mathematics, 
The Johns Hopkins University\\
Baltimore, MD 21218\\\smallskip
\\Department of Mathematics,
University of Oregon\\Eugene, OR 97403
\\\smallskip\\\rm Email:\stdspace\tt
ma2m@faraday.clas.Virginia.edu\quad jack@math.jhu.edu
\\sadofsky@math.uoregon.edu}

\asciiaddress{Department of Mathematics\\ 
Unversity of Virginia\\
Charlottesville, VA  22903\\
\\Department of Mathematics\\ 
The Johns Hopkins University\\
Baltimore, MD 21218\\
\\Department of Mathematics
University of Oregon\\Eugene, OR 97403
\\Email: ma2m@faraday.clas.Virginia.edu\\jack@math.jhu.edu
\\sadofsky@math.uoregon.edu}

\begin{abstract} 
We construct splittings of some completions of the $\mathbf{Z}/ (p)$--Tate
cohomology of $E (n)$ and some related spectra.  In particular, we
split (a completion of) $tE (n)$ as a (completion of) a wedge of $E (n-1)$'s as a
spectrum, where $t$ is shorthand for the fixed points of the
$\mathbf{Z}/ (p)$--Tate cohomology spectrum  (ie the Mahowald inverse
limit $\invlim{k}{(P_{-k} \wedge\Sigma  E (n))}$).  We also give a multiplicative
splitting of $tE (n)$ after a suitable base extension.
\end{abstract}

\asciiabstract{%
We construct splittings of some completions of the Z/(p)-Tate
cohomology of E(n) and some related spectra.  In particular, we
split (a completion of) tE(n) as a (completion of) a wedge of E(n-1)'s as a
spectrum, where t is shorthand for the fixed points of the
Z/(p)-Tate cohomology spectrum  (ie Mahowald's inverse
limit of P_{-k} smash SE(n) ).  We also give a multiplicative
splitting of tE(n) after a suitable base extension.}




\maketitlepage
%
%

\section{Introduction}\label{sec:Introduction}

\subsection{Preliminaries and notation}

We fix a prime $p$, and an integer $n$. 
We use $t$ to denote Mahowald's inverse limit
construction $tE = \invlim{k}{(P_{-k}\wedge \Sigma E)}$, where $P_{-k}$
stands for either $\mathbf{R}P^{\infty }_{-k}$ or its analogue for
$B\mathbf{Z}/ (p)$ when $p$ is odd; see \cite{MR:Root} or
\cite{Sad}.  This is an
abbreviation for the fixed points of the Greenlees--May Tate cohomology
functor; we write $tE$ for what would be denoted in
\cite{GrM} by $t_{\zp } (i_{*}E)^{\zp }$.  In particular 
if $E$ is a ring spectrum, then so is $tE$.

The starting point of this paper is a calculation on coefficients, which is
by now well known (see Lemma~\ref{t-le:Coefficients}).  If $E$ is a complex
oriented ring spectrum in which the series $[p] (x)$ is not a
zero divisor in $E^{*}\psb{x} = E^{*}\mathbf{C}P^{\infty }$, then  
\[
\pi _{-*}(tE) = E^{*}\lsb{x}/([p](x)),
\]
where $|x|=2$ and $E^{k} = \pi _{-k} (E)$. Here we use the awkward grading
to emphasize that the ring is related to $E^{*}\mathbf{C}P^{\infty }$
and $E^{*}B\mathbf{Z}/ (p)$.  Graded more conventionally, we have 
\[
\pi _{*} (tE) = E_{*}\lsb{x}/ ([p] (x))
\]
where $|x|=-2$.

Throughout this paper we
will use $R \lsb{x}$ to denote the ring of formal Laurent series over $R$ that
are allowed to be infinite series in $x$, but only finite in $x^{-1}$.  
If $E$ is a ring spectrum and $x$ an indeterminate in degree $-d<0$,
we write $E \lsb{x}$ for the ring spectrum given by $\invlim{i}{\bigvee
_{j\geq i}\Sigma ^{dj}E}$.   The multiplication is defined by the
inverse limit of the obvious maps
\begin{equation}\label{ring-wedge}
(\bigvee _{j\geq i_{0}}\Sigma ^{dj}E) \wedge (\bigvee _{k \geq
i_{1}}\Sigma ^{dk}E) \xrightarrow{} \bigvee _{l \geq i_{0}+i_{1}}\Sigma ^{ld}E.
\end{equation}
This gives $\pi _{*}[E \lsb{x}] = E_{*} \lsb{x}$. (Note that the
theoretical possibility of phantom maps in (\ref{ring-wedge}) means
the multiplicative structure there may not be unique, but we are
content to use whichever multiplicative structure might occur.
Actually work of Hovey and Strickland shows that there are no phantom
maps in this situation \cite{HoSt}, so the multiplication is defined
uniquely.) 
Note also that if $E$ is connective, 
\begin{equation}\label{limprod}
\invlim{i}{\bigvee
_{j\geq i}\Sigma ^{dj}E} = \prod _{j \in \ints }\Sigma ^{dj}E.
\end{equation}
In the succeeding, if we refer to a multiplication on $\prod _{j \in
\ints }\Sigma ^{dj}E$ it will be the one coming from the equivalence
(\ref{limprod}), and the ring structure on $\prod _{j \in \ints
}\Sigma ^{dj}R$ for a connective graded ring $R$ will be understood to
be the one given by the additive isomorphism with $R \lsb{x}$.

If $M$ is a flat module over $E_{*}$, we write $M \otimes E$ for the
spectrum representing the homology theory
\[
X \mapsto M \otimes _{E_{*}}E_{*} (X).
\]
This is a module spectrum over $E$ (again, it will not matter for us
that the module structure is well-defined
only up to phantom maps).  We also use $[p] (x)$ for the $p$--series of the formal
group law over $E$ given by the orientation.  When necessary we will
decorate $[p]$ with a subscript indicating the formal
group law.

We will work with a number of different spectra $E$, all
closely related.  The cohomology theory closest to $BP$ is $\bp{n}$,
the version of $BP$ with singularities satisfying $\bp{n}_{*} = \Zloc
[v_{1},\dots ,v_{n}]$ where $v_{i}$ is the $i$th Araki generator.
The Johnson--Wilson theory $\jw{n} = 
v_{n}^{-1}\bp{n}$ has the obvious coefficient ring obtained from
inverting $v_{n} \in \bp{n}_{*}$

We will also need to consider some simple variants of $\jw{n}$.  We list
the theories below and their coefficient rings.  They are all flat
over $\jw{n}$ and are thus determined by their coefficients:
\begin{align*}
\cjw{n} &= L_{K (n)}\jw{n} \qquad\text{so}\\
\cjw{n}_{*}& = \Zp \psb{v_{1},\dots,v_{n-1}}[v_{n},v_{n}^{-1}] 
            = \iadiccomp{(\jw{n}_{*})}{n-1}\\
(\gltjw{n})_{*} & = \cjw{n}_{*}[u]/ (u^{p^{n}-1}-v_{n})\\
&\text{\qquad\qquad\qquad and in later sections}\\
\ERT_{*} &= \jw{n}_{*}[\w]/ (\w ^{p^{n-1}-1}-v_{n-1}). 
\end{align*}
The ``$\wedge $'' is meant to suggest ``complete.'' 
The reader should note that our choice of $\gltjw{n}$ so that
$\gltjw{n,*} = \Zp \psb{u_{1},\dots ,u_{n-1}}[u,u^{-1}]$ conflicts with
the common convention that $\gltjw{n}$ denote the same theory extended
by the Witt vectors of $\mathbf{F}_{p^{n}}$.   We make no use of that
theory in this paper, so this should
cause no confusion.  The multiplicative structure on $\cjw{n}$ is given by the
composite
\begin{multline*}
L_{K (n)}\jw{n}\wedge L_{K (n)}\jw{n} \rightarrow L_{K (n)} (L_{K
(n)}\jw{n}\wedge L_{K (n)}\jw{n}) =\\
 L_{K (n)} (\jw{n}\wedge
\jw{n})\rightarrow L_{K (n)}\jw{n}.
\end{multline*}
The other two spectra are finite extensions of $\cjw{n}$, so the
multiplicative structure is determined in the obvious way from the
coefficients. 

\subsection{History}

Perhaps a few words of history are relevant here.  Lin's theorem
\cite{Lin,Gun} shows that for any finite CW--spectrum $X$, the map 
\begin{equation*}
X \rightarrow \invlim{k}{(P_{-k} \wedge \Sigma X)} = tX
\end{equation*}
is $p$--completion.   This is not true \cite{DM,DJKMW} for $E = BP$ or various
other $BP$--module spectra, though there is some predictable behavior, which
is quite different from that of finite complexes.  In particular, $tX$ is
quite large in these cases.

By a lovely but simple argument involving the thick subcategory
theorem, Mahowald and Shick \cite{MS:Root} showed (for $p=2$) that
if $X$ is a type $n$ finite complex, and $v$ is a $v_{n}$--self map, then $t 
(v^{-1}X) \simeq  *$.  This is the starting point for a series of
observations.   For example, it turns out that $tK (n) \simeq *$;
indeed \cite{GrS} proves that $t_{G}K (n) \simeq *$ for any finite
group $G$.  Hopkins has shown that if $X$ is of type $n$ then 
$t (L_{n}X) \simeq *$ if $X$ is type $n$ (this is proved in
\cite{H-S:Tate}).   Calculations of Hopkins  and Mahowald 
(including a proof when $n=1$) lead them to conjecture that if $X$ is
type $n-1$, then  
\begin{align}\label{Hopkins-Mahowald}
t (L_{n}X) = \pcomp{L_{n-1}X} \vee \Sigma ^{-1}\pcomp{L_{n-1}X}.
\end{align}
This is related to the chromatic splitting conjecture \cite{Hov}.

In light of (\ref{Hopkins-Mahowald}), it seems worthwhile to
investigate $tE (n)$.  In \cite{GrS} it is shown that $tE (n)$ is
$v_{n-1}$--periodic in an appropriate sense.  Using this, one
sees in \cite{H-S:Tate} that when $X$ is finite, the \emph{Bousfield class} 
of $t (L_{n}X)$ is compatible with the conjecture
(\ref{Hopkins-Mahowald}) (for $p=2$, this is also proved by different
techniques in \cite{CNL}).   The goal of the present paper is to give as
precise a description as possible of $tE (n)$ in terms of familiar
$v_{n-1}$--periodic spectra like $E (n-1)$.  We hope that with the
results of this paper in hand, it may be possible to make further
headway on (\ref{Hopkins-Mahowald}). 

\subsection{Results}

The first step, in section \ref{sec:Coefficients}, is a series of
calculations of $\pi _{*}tE$ in terms of more familiar objects for
various $E$.   It turns out to be useful to complete with
respect to $I_{n-1} = (p,\dots ,v_{n-2})$.  Since $tE (n)$ is already
$E (n-1)$ local, this amounts to localizing with respect to $K (n-1)$.  
This completion is in some sense not very dramatic; in particular if
$X$ is type $n-1$ as in (\ref{Hopkins-Mahowald}), the completion
leaves $t (L_{n}X)$ unchanged.   The main result, for the case of $\jw{n}$, is:

\begin{ITheorem}[Proposition \ref{En-coeff-iso}]
\label{it:En-coeffs-iso-equation}
There is an isomorphism of rings 
\begin{equation}
\iadiccomp{\bigl(t\jw{n}_{*}\bigr)}{n-1} =
\iadiccomp{\jw{n-1}_{*} \lsb{x}}{n-1}. 
\end{equation}
\end{ITheorem}

Section~\ref{sec:Module-spectrum}  gives
splittings of the spectra $v_{n}^{-1}t\bp{n}$ and $t\jw{n}$ in the manner
suggested by Theorem \ref{it:En-coeffs-iso-equation}.  For $t\jw{n}$
the result is:

\begin{ITheorem}[Theorem \ref{tate-splitting}] \label{it:tate-splitting}
There is a map of spectra
\begin{equation*}
\invlim{i\in \mathbf{N}}{\bigvee _{j \geq -i}\Sigma ^{2i}\jw{n-1}}
\rightarrow \iadiccomp{t\jw{n}}{n-1},
\end{equation*}
which after completion at $I_{n-1}$ (or
equivalently after localization with respect to $K (n-1)$) induces the
isomorphism of Theorem \ref{it:En-coeffs-iso-equation} on homotopy groups.
\end{ITheorem}

We emphasize that the map in Theorem \ref{it:tate-splitting} is not
multiplicative, even though the left hand side can be given a
multiplication as in equation (\ref{ring-wedge}) and the right hand
side is a ring spectrum by \cite{GrM},  and Theorem 
\ref{it:En-coeffs-iso-equation} gives an isomorphism of rings.  
Now the homotopy calculation in section \ref{sec:Coefficients}
also shows that the obvious formal group law over $t\jw{n}$ has height
$n-1$, so one might view the map 
\[
 \jw{n} \rightarrow t\jw{n}
\]
as a sort of Chern character, the classical Chern character being the
case $n=1$.  To make sense of this, one ought to use the
calculations of section \ref{sec:Coefficients} to construct a
map of ring spectra.  

To do this, we give an isomorphism between the natural formal
group law over $tE (n)$ and the Honda formal group law of
height $n-1$.  This can be done after a dramatic base extension,
but in section~\ref{sec:Isomorphism}  we show that it can then be done
in a canonical way.  Section~\ref{sec:Generalized-chern-character} uses 
section~\ref{sec:Isomorphism} and the Lubin--Tate theory of lifts to
construct an isomorphism of ring spectra between an extended version
of $tE (n)$ and an extended version of $E (n-1)$.

\begin{ITheorem}[Corollary \ref{iso-ring-spec}]\label{iso-ring-spec-theorem}
These is a canonical equivalence of complex oriented ring spectra
\[
\C_{\khuge}\hat{\tensor{\C_{\k\lsb{y}}}}TE \rightarrow
\C_{\khuge}\hat{\tensor{\Zp }} \gltjw{n-1}.
\]
\end{ITheorem}

The objects discussed in
Theorem~\ref{iso-ring-spec-theorem} are defined precisely in sections
\ref{sec:Isomorphism}  and \ref{sec:Generalized-chern-character}.
Briefly, the ring $\C_{\k\lsb{y}}$ is isomorphic to the $p$--adic
completion of 
$\Zp\lsb{\y}$ (see Lemma~\ref{t-le:CohenOfRes}).
The ring $\C_{\khuge }$
is a complete discrete valuation ring with maximal ideal $(p)$ and
residue field $\khuge $.  $TE$ is a finite extension of $tE (n)$
completed at $I_{n-1}$ and defined in equation (\ref{bigtate}), and
$\gltjw{n-1}$ is defined as usual to be a completed $2$--periodic
extension of $\jw{n-1}$ as 
usual so that (\ref{def:gltjw}) holds.

The authors thank Johns Hopkins University, at which the conversations
leading to this paper took place, and Max-Planck-Institut f\"{u}r
Mathematik, at which the last fragments were committed to paper.
The authors were partially supported by the NSF.

\section{Calculations on homotopy groups}
\label{sec:Coefficients}

\subsection{The homotopy of $tE$ for $E$ complex oriented}

We assume that $E$ is complex oriented, and that $[p] (x)$ is not a
$0$ divisor in $E^{*}\psb{x}$.  In \cite[\S 16]{GrM} it is shown 
that $tE$ 
is an inverse limit of $E$ smashed with Thom spectra over
$B\mathbf{Z}/ (p)$.  In particular, \cite[16.3]{GrM} shows  that $\pi
_{-*}tE = E^{*} \lsb{x}/ ([p] (x))$ as a 
module over $E^{*}\psb{x}/ ([p] (x)) = E^{*} (B\mathbf{Z}/ (p)_{+})$.

Now \cite[3.5]{GrM} proves that the map 
\[
E^{*} (B\mathbf{Z}/ (p)_{+}) \rightarrow \pi _{-*}tE
\]
is a map of rings.  It follows that the element ``$x^{-1}$'' in $\pi
_{-*}tE = E^{*} \lsb{x}/ ([p] (x))$ satisfies $x \cdot x^{-1} = 1$ in
the ring $\pi _{*}tE$, so really is the inverse of $x$.  From this one
concludes that 

\begin{Lemma} \label{t-le:Coefficients}
There is an isomorphism of rings
\begin{equation}\label{coefficient-calculation} 
\pi _{-*}(tE) = E^{*}\lsb{x}/([p](x)).
\end{equation}
\end{Lemma}

\subsection{An isomorphism after completion}

Recall $R \lsb{x}$ is $S^{-1}R\psb{x}$ where $S$ is the
multiplicatively closed subset generated by $x$.  If $R$ is
graded connected, and $x$ has degree $-2$, then additively $R\lsb{x} =
\prod_{k\in \ints }\Sigma ^{2k}R$.   As remarked in the introduction,
wherever this notation occurs, we 
will also consider $\prod_{k\in \ints }\Sigma ^{2k}R$ to have the ring
structure given by this isomorphism.  Similarly, when $E$ is a
connective ring spectrum, we will consider the ring structure on
$\prod_{k\in \ints }\Sigma ^{2k}E = \invlim{i}{\bigvee
_{k\geq i}\Sigma ^{2k}E}$ as in (\ref{ring-wedge}).

The following observation motivates
Conjecture~1.6 of \cite{DM} (which was corrected as Conjecture~1.2
of \cite{DJKMW}).

\begin{Proposition} \label{t-pr:bp<n>-coeff-iso}
\begin{eqnarray}\label{bp<n>-iso}
\pi _{*}(t\bp{n}) & = & \bp{n}_{*}\lsb{x}/([p](x)) \nonumber \\
 & \simeq & \prod_{k \in \ints } \Sigma
^{2k}\pcomp{\bp{n-1} }_{*}\\
 & \simeq & \pcomp{\bp{n-1}}_{*}\lsb{x}. \nonumber  
\end{eqnarray}
\end{Proposition}
\begin{proof}
The first line is equation (\ref{coefficient-calculation}).
We wish to simplify the right hand side.

 For convenience, we introduce degree $0$ elements
\[ 
  w_{i} = v_{i}x^{p^{i}-1}.
 \]
Then the ideal in the quotient above is generated by the relation
\[
px +_{F} w_{1}x +_{F} w_{2}x +_{F}\dots +_{F} w_{n}x = 0
\]
where $F$ is the formal group law over $\bp{n}_{*}$ induced by the
orientation from $BP$.  Expanding the formal sums and dividing by $x$,
we can write the relation in the form
\[
    w_{n} = - (p + \mbox{ a formal series in } w_{1},...,w_{n}),
\] 
and this equation implies
\begin{equation}\label{wn-image} 
     w_{n} \equiv  - (p + w_{1} + ... + w_{n-1}) \pmod{\mbox{decomposables}}. 
 \end{equation}
By iterating this relation some finite number of times (depending on
the multi-index $(i_{0},\dots ,i_{n-1})$), we produce a polynomial
$W_{n}$ so that 
\[
     w_{n} = W_{n}(w_{1},...,w_{n-1}) \in {\bf
Z}/ (p^{i_{0}})[w_{1},\dots ,w_{n-1}]/ (w_{1}^{i_{1}},\dots ,w_{n-1}^{i_{n-1}}).
\]
We thus produce a power series $W_{n} \in {\bf Z}_{p}[[w_{1},\dots
,w_{n-1}]]$ so that 
\[ 
     w_{n} = W_{n}(w_{1},...,w_{n-1}) \in {\bf Z}_{p}[[w_{1},\dots
,w_{n-1}]] \subseteq {\bf
Z}_{p}[v_{1},\ldots,v_{n-1}]\psb{x}.
 \]
Thus the map 
\[ 
     \pi _{*}(t\bp{n}) \rightarrow
\pcomp{\bp{n-1}}_{*}\lsb{x}
 \]
which sends $v_{i}$ to $v_{i}$ ($i<n$), $x$ to $x$ and $v_{n}$ to
$x^{-(p^{n}-1)}W_{n}(w_{1},...,w_{n-1})$ is a 
well defined ring map.  There is a map 
\[
\bp{n-1}_{*}\lsb{x} \rightarrow \pi _{*}(t\bp{n})
\]
defined in the obvious way, that becomes an inverse map upon extending
it to $\pcomp{\bp{n-1}}_{*}\lsb{x}$.  This extension exists, because
the relation
\[
[p] (x) = 0
\]
together with Araki's formula for the $p$--series
\cite[A2.2.4]{Ravenel:GreenBook} allows us to write 
\begin{equation}\label{formula-for-p}
px = [-1] (v_{1}x^{p}+_{F}\dots +_{F}v_{n}x^{p^{n}})
\end{equation}
where $F$ is the formal group law on $\bp{n}$ induced by the complex
orientation from $BP$.
Now dividing (\ref{formula-for-p}) by $x$ lets us to write any
$p$--adic integer as a power series in $x$.
\end{proof}

The conjecture of \cite{DM,DJKMW} is that the isomorphism
\eqref{bp<n>-iso} is the effect on homotopy groups of an equivalence
of spectra 
\[ 
t\bp{n} \simeq \prod_{k \in {\bf Z}}
\Sigma ^{2k} \pcomp{\bp{n-1}}.
 \]

This is proved for $n=1$ in \cite{DM} and for $n=2$ in \cite{DJKMW}
using the Adams spectral sequence.
Of course in general an isomorphism in homotopy doesn't give an
isomorphism of spectra, but it does if it is an isomorphism of
$MU_{*}$--modules satisfying Landweber exactness.  This suggests
inverting $v_{n}$ to produce an isomorphism of spectra.  We shall show
that there is an isomorphism of the expected form after an appropriate
completion, taking as our starting point the spectra $v_{n}^{-1}
(t\bp{n})$ and $tE (n)$.  The first point is that these spectra are
not the same.

The map $\bp{n} \rightarrow \jw{n}$ gives a map
$t\bp{n} \rightarrow t\jw{n}$, and since $v_{n}$ is a unit on the
right, this gives a map $v^{-1}_{n}t\bp{n} \rightarrow t\jw{n}$.  On
homotopy we get
\begin{eqnarray*}
\pi _{*} ( v_{n}^{-1}(t\bp{n})) & = &
v_{n}^{-1}[\bp{n}_{*}\lsb{x}/([p](x))]\\
 & \hookrightarrow & \jw{n}_{*}\lsb{x}/([p](x))\\
 & = &  (t\jw{n})_{*} = \pi _{*}t (v_{n}^{-1}BP\langle n  \rangle).
\end{eqnarray*}
To see that the inclusion is proper, set $r = |v_{n}|/|v_{1}| + 1$
and notice that the series
\[ 1 + v_{n}^{-1}v_{1}^{r}x^{p-1} + v_{n}^{-2}v_{1}^{2r}x^{2(p-1)} +
\ldots  \]
is an element of $\jw{n}_{*}\lsb{x}/([p](x))$ but not of $v_{n}^{-1}\bp{n}_{*}\lsb{x}/([p](x))$.
This reflects the fact that $t$ need not commute with direct limits
(nor in fact does it generally commute with inverse limits).

We treat the case $v_{n}^{-1} (t\bp{n})$ first.   


We would like to extend the isomorphism of (\ref{bp<n>-iso}) to an
isomorphism 
\begin{equation}
\label{vbp<n>-iso}
\pi_{*}v_{n}^{-1}t\bp{n} \rightarrow
v_{n-1}^{-1} \pcomp{\bp{n-1}}_{*}\lsb{x}.
\end{equation}
For the map to exist we need the image of $w_{n}$ under the
isomorphism in (\ref{bp<n>-iso}) to be a unit after
inverting $v_{n-1}$ on the range of that map.  This is false for $n >
1$.  We can see 
this by checking that the  image of $w_{n}$ after inverting $v_{n-1}$
is not a unit even modulo
$(v_{1},\ldots,v_{n-2})$.  The range of
(\ref{vbp<n>-iso}) in dimension $0$, modulo $(v_{1},\ldots,v_{n-2})$, is $\Zp \psb{w_{n-1}}[w_{n-1}^{-1}]$.
We have $\im (xw_{n}) = - (px+_{F}w_{n-1}x)$, so 
\[
\im (w_{n}) = - (p+w_{n-1}) \pmod{pw_{n-1}\Zp \psb{w_{n-1}}}. 
\]
Since $w_{n-1}$ is a unit in the range, the image of $w_{n}$ is a unit
if and only if the image of $w_{n}$ divided by $w_{n-1}$, which we will
abuse notation by writing $w_{n}/w_{n-1}$, is a unit.  
Continuing to work modulo $(v_{1},\dots ,v_{n-2})$,
\[
w_{n}/w_{n-1} = - (1+\frac{p}{w_{n-1}})  \pmod{p\Zp \psb{w_{n-1}}}. 
\]
Now if we examine the inclusion of rings
\[
\eta \co \Zp \psb{w_{n-1}}[w_{n-1}^{-1}] \subseteq \pcomp{(\Zp \psb{w_{n-1}}[w_{n-1}^{-1}])},
\]
we see that since $w_{n}/w_{n-1} = -1+p\varepsilon $, $w_{n}/w_{n-1}$
is a unit in the completion, with inverse given by $- (1+p\varepsilon +p^{2}\varepsilon ^{2}+\dots)$.
If we write this element as a power series in $w_{n-1}$ plus a power
series in $w_{n-1}^{-1}$, we see that $(w_{n}/w_{n-1})^{-1}$ is the sum of
a power series in $w_{n-1}$ with 
\[
\frac{a_{1}}{w_{n-1}}+\frac{a_{2}}{w_{n-1}^{2}}+\frac{a_{3}}{w_{n-1}^{3}}+\dots 
\]
where $a_{i} = (-1)^{i}p^{i} \pmod{p^{i+1}}$.
The power series in $w_{n-1}$ is  neccessarily in the image of $\eta
$, but
$\frac{a_{1}}{w_{n-1}}+\frac{a_{2}}{w_{n-1}^{2}}+\frac{a_{3}}{w_{n-1}^{3}}+\dots$
is clearly not.  It follows that $(w_{n}/w_{n-1})^{-1}$ is not in the
image of $\eta $, and hence $w_{n}/w_{n-1}$ is not a unit in $\Zp \psb{w_{n-1}}[w_{n-1}^{-1}]$.

Similarly, if the map existed, it could be an isomorphism only if $v_{n-1}$ 
were a unit in the domain after inverting $v_{n}$.  This also is not the
case for $n > 1$, for similar reasons. When $n=1$ both of these conditions 
are met, the map above exists, and is an isomorphism.

\begin{Proposition}\label{funny-iso}
If the map from (\ref{bp<n>-iso}) is completed at the ideal 
\[ 
I_{n-1}
= (p,v_{1},\ldots,v_{n-2})
 \]
then there is an isomorphism
\begin{equation}
\label{comp-iso}
\iadiccomp{(\pi_{*} v_{n}^{-1}t\bp{n})}{n-1} \rightarrow
  \iadiccomp{( v_{n-1}^{-1} [\pcomp{\bp{n-1}}_{*}\lsb{x}])}{n-1}.
\end{equation}
\end{Proposition}

\begin{proof} We follow the isomorphism of (\ref{bp<n>-iso}) 
\begin{eqnarray*}
\pi _{*}(t\bp{n}) & = & \bp{n}_{*}\lsb{x}/([p](x)) \\   
                  & \simeq & \pcomp{\bp{n-1}}_{*}\lsb{x}
\end{eqnarray*}
by the inclusion 
\[
  \pcomp{\bp{n-1}}_{*}\lsb{x} \rightarrow 
  \iadiccomp{(v_{n-1}^{-1} [\pcomp{\bp{n-1}}_{*}\lsb{x}])}{n-1}
\]
into the module obtained by inverting $v_{n-1}$ and then completing.  Now note
that $w_{n}$ has image as in (\ref{wn-image}).  So dividing by
$v_{n-1}x^{p^{n-1}}=w_{n-1}$ we see the image of $w_{n}$ divided by
$w_{n-1}$ is  $-1$ plus terms in the ideal $I_{n-1}$ and terms in
$x\pcomp{\bp{n-1}}_{*}\psb{x}$. This is a unit, so since we have
inverted $v_{n-1}$, the image of $w_{n}$ is a unit.  This allows us to
extend our map to the domain given by inverting $v_{n}$.  Since
the range is complete, we can also extend to the completion.

On the other hand, a similar argument allows us construct the inverse
map from the inverse map of (\ref{bp<n>-iso}).
\end{proof}

Now we shall show that the isomorphism in
(\ref{comp-iso}) is induced by a map of spectra
\[
\iadiccomp{(v_{n}^{-1}t\bp{n})}{n-1}  \xrightarrow  {}
\iadiccomp{( v_{n-1}^{-1}\prod_{\ints }\Sigma ^{2k}\bp{n-1})}{n-1}.
\]
%
%
%
Both sides are ring spectra with obvious $MU$--module structures.  
They would be isomorphic as $MU$--algebras by the
Landweber exact functor theorem if we could make them $MU$--algebras 
so that the coefficient isomorphism preserves the map from $MU_{*}$.  

Let $R$ be the
ring spectrum on the right.  In order to construct a map of ring
spectra inducing the isomorphism in (\ref{comp-iso}), it is
necessary that the FGL induced by 
\[ 
v_{n}^{-1}t\bp{n}_{*} \rightarrow
R_{*} \rightarrow R_{*}/(p,\ldots,v_{n-2})
 \]
be isomorphic to the ``usual'' FGL on $R_{*}/(I_{n-1})$
(induced by $\bp{n-1}_{*} \rightarrow R_{*}$).
We cannot demonstrate such an isomorphism, but we can exhibit an
isomorphism of spectra (that preserves neither the $MU$--module
structure, nor the multiplicative structure).   We do this in
section~\ref{sec:Module-spectrum}. 


The situation for $t\jw{n}$ is very similar.  We
can attempt to construct a map 
\begin{equation}\label{homotopy-splitting} 
t\jw{n}_{*} \rightarrow \pcomp{\jw{n-1}}_{*}\lsb{x}
 \end{equation}
as in (\ref{bp<n>-iso}), but we immediately
run into the problem that $w_{n}$, and hence $v_{n}$, does not go to
a unit for $n > 1$.  Also as before, $v_{n-1}$ is not a unit on the
left. The solution is the same; after completing both sides at $I_{n-1}$ we
can construct an isomorphism.  
\begin{Proposition}\label{En-coeff-iso}
The map of (\ref{bp<n>-iso}) extends to an isomorphism
\begin{equation}\label{En-coeffs-iso-equation}
\iadiccomp{\jw{n}_{*} \lsb{x}/[p] (x)}{n-1} =
\iadiccomp{\jw{n-1}_{*} \lsb{x}}{n-1} 
\end{equation}
where 
\[
\iadiccomp{\pi _{*} (tv_{n}^{-1}\bp{n})}{n-1} = \iadiccomp{\pi
_{*}t\jw{n}}{n-1} = \iadiccomp{\jw{n}_{*} \lsb{x}/[p] (x)}{n-1}
\]
and
\[
\iadiccomp{\jw{n-1}_{*} \lsb{x}}{n-1} =
\iadiccomp{(v_{n-1}^{-1}\bp{n-1}) \lsb{x})}{n-1}.
\]
\end{Proposition}

In the next section, we use this calculation to construct an
isomorphism of spectra.  We are able to do
this without showing the corresponding formal group laws are
isomorphic, so we do not get an isomorphism of $MU$--modules.  In section
\ref{sec:Generalized-chern-character}, we extend scalars suitably to
construct an isomorphism of formal group laws, yielding isomorphisms
of $MU$--algebra spectra.

%

\section{Structure of $\EFP{n} $ as an $\EFP{n-1} $--module
spectrum} \label{sec:Module-spectrum}

We begin with an algebraic observation.

\begin{Lemma}\label{algebraic-splitting}
Let $(A,\mathfrak{m})$ be a complete local ring, $k = A/\mathfrak{m}$.  Let $M$ be an
$A$--module such that the map $\mathfrak{m} \otimes x \mapsto \mathfrak{m}x$ induces
isomorphisms 
\[
(\mathfrak{m}^{r}/\mathfrak{m}^{r+1})\otimes _{k} (M/\mathfrak{m}M) = \mathfrak{m}^{r}M/\mathfrak{m}^{r+1}M.
\] 
Let $I$ be an index set for a vector space basis of $M/\mathfrak{m}M$,
so that $M/\mathfrak{m}M \simeq \bigoplus _{I}k$.  Then
there is a map 
\[
\bigoplus _{I}A \xrightarrow{\varepsilon } M
\]
which is an isomorphism when completed at $\mathfrak{m}$.
\end{Lemma}
\begin{proof}
We construct a map as follows.  Let $\{\bar{x}_{\alpha } \}, \alpha \in
I$ be a basis for $M/\mathfrak{m}M$ as a $A/\mathfrak{m}$--vector space.
Let $\{x_{\alpha } \}$ be lifts to $M$.  Then we map $\oplus _{I}A$
to $M$ by $\varepsilon (1_{\alpha }) = x_{\alpha }$.

Our hypotheses imply that this map gives an isomorphism on the
associated graded with respect to the filtration induced by powers of
$\mathfrak{m}$.  This implies that the completion of the map is an isomorphism.
\end{proof}

Note that if the index set $I$ in Lemma~\ref{algebraic-splitting} is
infinite, $M$ 
will not generally be isomorphic to a free $A$--module.  For
example if $A = \mathbf{Z}_{p}$ and $M =
\pcomp{(\oplus_{\mathbf{N}}\mathbf{Z}_{p})}$  then $M$ is not a free  
$\mathbf{Z}_{p}$--module.  To see this, observe that no free
$\mathbf{Z}_{p}$--module of infinite rank can be $p$--complete.  

Next we recall the result proved in \cite[Theorem 4.1]{H-S:Pic}.  The
identification $BP_{*}BP = BP_{*}[t_{1},t_{2},\dots ]$ gives a splitting
\begin{equation}\label{bp-splitting}
\psi \co BP \wedge BP \simeq \bigvee _{R}\Sigma ^{|R|}BP
\end{equation}
where $R$ ranges over multi-indices of non-negative integers (with
only finitely many positive coordinates) $R= (r_{1},r_{2},r_{3},\dots
)$, $t^{R} = t_{1}^{r_{1}}t_{2}^{r_{2}}\dots $  and 
\[
|R| = |t^{R}| = 2 (r_{1} (p-1)+r_{2} (p^{2}-1)+r_{3} (p^{3}-1)+\dots).
\]
To build a map  from right to left of (\ref{bp-splitting}), take $\Sigma ^{|R|}BP$
to $BP \wedge BP$ by using the homotopy class $t^{R} \in \pi _{|R|} (
BP\wedge BP)$, smashing on the left with $BP$ and then multiplying the
left pair of $BP$'s:
\begin{equation}\label{summand-inclusion}
BP \wedge S^{|R|} \xrightarrow{1_{BP} \wedge t^{R}} BP \wedge BP \wedge BP \xrightarrow {\mu
\wedge 1_{BP}} BP \wedge BP.
\end{equation}
The wedge over all $R$ of these maps is an isomorphism on homotopy
groups, so is invertible, and $\psi $ is that inverse.

Theorem 4.1 of \cite{H-S:Pic} states that the composite 
\begin{equation}\label{long-map}
BP \xrightarrow {\eta _{R} } BP \wedge BP \xrightarrow{\psi } \bigvee
_{R}\Sigma ^{|R|}BP \xrightarrow{\rho } \bigvee _{R \in \mathcal{R}}\Sigma
^{|\sigma R|}BP \xrightarrow{\theta } \bigvee _{R \in \mathcal{R}}\Sigma
^{|\sigma R|}BP\langle j  \rangle
\end{equation}
is a homotopy equivalence after smashing with a type $j$ spectrum and
inverting $v_{j}$.  The map $\rho $ is induced by leaving out wedge summands,
and $\theta $ by the usual reduction of ring spectra $BP \rightarrow BP\langle j  \rangle$.
We use the notation of \cite{H-S:Pic}, derived from
\cite{JW1}:  $\mathcal{R}$ is the set of multi-indices with the first
$j-1$ indices $0$, and $\sigma R = (p^{j}e_{j},p^{j}e_{j+1},\dots
)$.   
%
%

We shall use 
the following facts about $MU$--module spectra $M$:
\begin{equation}\label{vn-limit}
M \wedge L_{j}Z = v_{j}^{-1} M \wedge Z =  M \wedge v_{j}^{-1}Z
\end{equation}
when $Z$ is a finite type $j$ spectrum (which follows from \cite[Theorem 1]{Rav:Geo});
\begin{equation}\label{iadic-completion}
\iadiccomp{M}{j-1} = \invlim{Z}{M \wedge Z}
\end{equation}
where $Z$ runs through finite type $j$
spectra under $S^{0}$.  Use of the nilpotence theorem, \cite{DHS},  is
required to produce 
a sufficient supply of such $Z$ and to ensure there are enough maps
between them, as in \cite[Proposition~3.7]{MSad}.  For our purposes
equation (\ref{iadic-completion}) can be taken to be a definition, but
see the remark below.  Inverting $v_{j}$, we get 
\begin{equation*}
L_{K (j)}M = \invlim{Z}{(M \wedge v_{j}^{-1}Z)} = \iadiccomp{(v_{j}^{-1} M)}{j-1} .
\end{equation*}
This last equality is by equation (\ref{vn-limit}) combined with, say,
the proof of \cite[Proposition 7.4]{HMS} which verifies that $L_{K
(n)}X = \invlim{Z}{(X \wedge L_{n}Z)}$ as $Z$ runs over finite type
$n$ spectra under $S^{0}$.

So if $v_{j}$ is already a unit on $M$ then
\begin{equation*}
L_{K (j)}M = \iadiccomp{M}{j-1}
\end{equation*}
and we will generally use the first notation rather than the second below.

\rk{Remark}
Although we use (\ref{iadic-completion}) as the
definition of $\iadiccomp{M}{j-1}$, there are other approaches that
can be taken for specific $M$.  For certain $MU$--modules $M$, one can
define a spectrum $\tilde{M}$ with homotopy groups $\iadiccomp{(\pi
_{*}M)}{j-1}$ by 
using Landweber exactness.  In those cases, also by using Landweber 
exactness, one can prove that $\tilde{M} \simeq \iadiccomp{M}{j-1}$ as
given in (\ref{iadic-completion}).
If $M$ has enough structure, one may also be able to define a
completion of $M$ using either Baas--Sullivan bordism with
singularities \cite{Baas}, or structured ring spectra using the
techniques of \cite{MR97h:55006}.  In either case, one can use the nilpotence
theorem to verify that the construction is homotopy equivalent to the
one in (\ref{iadic-completion}).

Each map in (\ref{long-map}) except $\eta _{R}$ is a map of left
$BP$--modules.  Recall $I_{j}$ is the ideal $(p,\dots ,v_{j-1})
\subseteq BP_{*}$.  Since $I_{j}$ is invariant, $\eta _{R}
(I_{j}) = \eta _{L} (I_{j})$, and thus each map in
(\ref{long-map}) is compatible on homotopy with the $I_{j}$--adic
filtration. Let $Z$ be a finite type $j$ spectrum; then smashing
(\ref{long-map}) with $L_{j}Z$ gives an equivalence, and thus an
$E(j)$--module structure on $BP \wedge
v_{j}^{-1}Z = BP \wedge L_{j}Z$ 
or, by taking inverse limits, a (possibly not unique) module structure
on $L_{K (j)}BP = \invlim{Z}{BP \wedge v_{j}^{-1}Z}$.

We prove the following proposition as a warm-up to our additive
splitting of the Tate cohomology spectrum.  The construction is
a general method for splitting $L_{K (j)}F$ when $F$ is a
nice enough $BP$--module, given the splitting of $L_{K (j)}BP$.
A theorem equivalent to Proposition~\ref{Ej-splitting} is proved as
\cite[Theorem 4.7]{H-S:Pic}.

\begin{Proposition}\label{Ej-splitting}
If $n > j$, there is a map 
\begin{equation}\label{Ej-splitting-equation}
s\co \bigvee _{V}\Sigma ^{|V|} E (j) \xrightarrow{}  L_{K (j)}E (n)
\end{equation}
which gives an equivalence after completing with respect to $(p,\dots
,v_{j-1})$, or equivalently, after localizing with respect to $K (j)$.  
The index $V$ runs through the monomials in 
\[
\mathbf{F}_{p}[v_{j+1},\dots
,v_{n-1},v_{n}^{\pm 1}].
\]
\end{Proposition}
\begin{proof}
$L_{K (j)}E (n)$ is an $L_{K (j)}BP$--module, and thus an $E
(j)$--module.  Let 
$\mathfrak{m} = I_{j} = (p,\dots ,v_{j-1})$. Since
\begin{align}\label{}
L_{K (j)}E (n) = \invlim{Z}{(E (n) \wedge v_{j}^{-1}Z)}
\end{align}
where $Z$ runs through finite type $j$ spectra,
\[
\pi _{*}L_{K (j)}E (n) = \mathbf{Z}_{p}\psb{v_{1},\dots
,v_{j-1}}[v_{j}^{\pm 1},v_{j+1},\dots ,v_{n-1},v_{n}^{\pm 1}].
\]

Now since $\mathfrak{m}$ is an invariant ideal 
$\mathfrak{m}\pi _{*}L_{K (j)}BP$ is well-defined, whether we
think of $\mathfrak{m}$ as an ideal in $E (j)_{*}$ acting on $M=\pi _{*}L_{K (j)}E (n)$ via the $E (j)$--module
action on $L_{K (j)}BP$ or as an ideal in $BP_{*}$ acting via the associated
localization map to $\pi _{*}L_{K (j)}BP$.

We calculate that the associated graded to the $\mathfrak{m}$--adic
filtration on $M$ is 
\begin{equation*}
E_{0}M = \oplus_{V} \Sigma ^{|V|}\mathbf{F}_{p}\psb{v_{0},\dots ,v_{j-1}}[v_{j}^{\pm 1}]
\end{equation*}
where $V$ runs through monomials in $\mathbf{F}_{p}[v_{j+1},\dots
,v_{n-1},v_{n}^{\pm 1}]$.

$M$ satisfies the hypotheses of Lemma~\ref{algebraic-splitting} (for
the left $\mathfrak{m}$--structure which comes from the map 
$BP \rightarrow E (n) \rightarrow L_{K (j)}E (n)$), where 
\[
A = \pi
_{*}\cjw{j} = \mathbf{Z}_{p}\psb{v_{1},\dots ,v_{j-1}}[v_{j}^{\pm 1}].
\]
 
Now the splitting of $L_{K (j)}BP$ gives $L_{K (j)}E (n)$  an $\cjw{j}$--structure,
and thus an associated $A$--action; the
$\mathfrak{m}$--adic filtration is the same, so $M$ also satisfies 
the hypotheses of Lemma~\ref{algebraic-splitting}  for that 
$\mathfrak{m}$--adic filtration. 

We can now make the usual homotopy theoretic argument: take generators of
$M/\mathfrak{m}M$ and lift them to elements of $M$.  Use the $E (j)$
structure of $L_{K (j)}E (n)$ to make maps $E (j) \rightarrow L_{K
(j)}E (n)$ realizing these generators on the unit of the ring spectrum
$E (j)$.  This gives a map 
\[
\bigvee_{V} \Sigma ^{|V|} E (j) \rightarrow L_{K (j)}E (n).
\]
By Lemma~\ref{algebraic-splitting}, this map induces an isomorphism on
homotopy groups after completion with respect to $(p,\dots ,v_{j-1})$,
ie after applying $L_{K (j)}$, which leaves the right hand side unchanged.
\end{proof}

We apply techniques similar to those used in
Proposition~\ref{Ej-splitting}  to prove the following.

\begin{Theorem}\label{tate-splitting}
There is a map of spectra
\begin{equation*}
\invlim{i\in \mathbf{N}}{\bigvee _{j \leq i}\Sigma ^{2i}E (n-1)} \rightarrow \iadiccomp{tE (n)}{n-1}
\end{equation*}
that becomes an isomorphism on homotopy groups after completion at $I_{n-1}$ \textup{(}or
equivalently after localization with respect to $K (n-1)$\textup{)}.
\end{Theorem}
\begin{proof}
We proceed as in the proof of Proposition~\ref{Ej-splitting}.  We have
given $M=\pi _{*}\iadiccomp{tE (n)}{n-1}$ as a $BP_{*}$--module in
Proposition~\ref{En-coeff-iso}.  It satisfies the hypotheses of
Lemma~\ref{algebraic-splitting} with 
$\mathfrak{m} = (p,\dots ,v_{n-2}) \subseteq BP_{*}$, $A =
\cjw{n-1}_{*}$.  Since
$\iadiccomp{tE (n)}{n-1} = L_{K (n-1)}tE (n)$, we have an $L_{K(n-1)}BP$--action and hence an $E (n-1)$--action.  

As above, the two available $\mathfrak{m}$--adic filtrations on 
$\pi
_{*}\iadiccomp{t E (n)}{n-1}$ (one from 
$BP_{*}$, the other from $E (n-1)_{*}$ via the $E(n-1)$--structure on
$L_{K (n-1)}BP$) are the same since $\mathfrak{m}$ 
is an invariant ideal.  We now proceed in a slightly different manner.
Note that  
\begin{equation}\label{vector-space}
M/\mathfrak{m} = K (n-1)_{*} \lsb{x} = \bigoplus _{i \in I}\Sigma ^{c_{i}}K (n-1)_{*}
\end{equation}
for some indexing set $I$, since $K (n-1)_{*}$ is a graded field.  We could apply
Lemma~\ref{algebraic-splitting}  and proceed as before, but we would
actually like better control over our expression for $\iadiccomp{tE
(n)}{n-1}$.  In particular, the index set $I$ in equation
(\ref{vector-space}) must be uncountable, but we would like to find a
countable set of topological generators for $\iadiccomp{tE
(n)}{n-1}$.  In fact, we would like these
generators to correspond to the (positive and negative) powers of
$x$. 

To accomplish this, we first recall \cite[Theorem 16.1]{GrM} which
states that $tE = \invlim{i}{[(B\mathbf{Z}/(p))^{- i\xi }\wedge \Sigma
E]}$, where $\xi $ is the usual complex line bundle.  
Recall also that the Thom class of $(B\mathbf{Z}/(p))^{- i\xi}$ is in
dimension $-2i$, and is not torsion.  In fact the spectrum
$(B\mathbf{Z}/(p))^{- i\xi}$ has a CW--structure with exactly one cell
in each dimension greater than or equal to $-2i$. Since the $-2i$--cell
generates a non-torsion element of homology, the attaching
map of the $-2i+1$ cell to the $-2i$ skeleton is null.  So the cell in
dimension $-2i+1$ is spherical, and the
inclusion of that $-2i+1$ cell, smashed with the unit of $E$ when $E$
is a ring spectrum, gives 
\[
x^{- i+1} \in \pi _{-2i+2}[(B\mathbf{Z}/(p))^{- i\xi} \wedge \Sigma E].
\]
%

Now, we take $x^{j} \in \pi
_{2j}\iadiccomp{t\jw{n}}{n-1}$, and use the $\cjw{n-1}$--structure to
construct a sequence of maps 
\begin{equation}\label{partial-wedge-mapped-in}
\bigvee _{j \geq -i}\Sigma ^{2j}\cjw{n-1} \xrightarrow{} \iadiccomp{t\jw{n}}{n-1}.
\end{equation}
We make a map $\mu _{-i}$ by composing the map of
(\ref{partial-wedge-mapped-in}) with the map 
\begin{equation*}
\iadiccomp{t\jw{n}}{n-1} \rightarrow \iadiccomp{(B\mathbf{Z}/
(p))^{- (i+1)\xi }\wedge \Sigma E}{n-1} 
\end{equation*}
given by \cite[Theorem 16.1]{GrM}.

Taking inverse limits of the maps $\mu _{-i}$ gives a map 
\begin{equation*}
\invlim{i}{(\bigvee _{j \leq i}\Sigma ^{2j}\cjw{n-1})}
\xrightarrow{f} \iadiccomp{t\jw{n}}{n-1}.
\end{equation*}
This map defines an isomorphism on the associated graded modules with respect to
$\mathfrak{m}$.   It follows that $f$ is an equivalence after
completion, that is 
\begin{equation*}
\iadiccomp{\bigl[\invlim{i}{(\bigvee _{j \leq i}\Sigma ^{2j}\cjw{n-1})}\bigr]}{n-1}
= \iadiccomp{t\jw{n}}{n-1}.
\end{equation*}
\vglue -1truecm
\end{proof}

By a very similar argument one can prove the analog to
Proposition~\ref{funny-iso}: 
\begin{Proposition}\label{funny-equiv}
There is an equivalence of spectra
\begin{equation*}
[ v_{n-1}^{-1} \invlim{i}{(\bigvee _{j\leq i}\Sigma ^{2j}\bp{n-1})}]\Iwadic
 \rightarrow (v_{n}^{-1} t\bp{n})\Iwadic .
\end{equation*}
\end{Proposition}
We leave the proof to the interested reader.  

%
%

Given all the completions that occur in this section and in section~\ref{sec:Coefficients}, one might hope that by using some other,
already complete theory like $\cjw{n}$ or $\ltjw{n}$, we could prove a
theorem with a simpler statement.  This is unfortunately not the case.  
There are similar results for these theories,
but even if $E$ is complete with respect to $(p,v_{1},\dots
,v_{n-2})$, $tE$ will generally not be, and so will need to be completed
again.  There are variants of Theorem~\ref{tate-splitting} for these
other spectra as well, but the statement is not simpler.

\section{A Honda coordinate on the formal group over $\Tate E
$} \label{sec:Isomorphism}

In this section we shall take $p$ to be an odd prime and $n>1$ to be an
integer.  It will ease the superabundance of superscripts to use the
abbreviation $q=p^{n-1}$.

We define a number of formal group laws in this section that are used
in the remainder of the paper.  For reference, we list them here.

\begin{enumerate}
\item [] $G (n)$: a homogeneous formal group law of degree $2$ on $E
(n)_{*}$ induced by the usual orientation of $E (n)$.
\item [] $\GLT{n}$: a twist of the pushforward of $G (n)$ to $\pi_{*}\gltjw{n}$ 
by the element $u$.  This is homogeneous of degree $0$.
\item [] $\GNice$: a twist of the pushforward of $G (n)$ to $\ERT _{*}$ of
$G (n)$ by $w$.  This is also homogeneous of degree $0$.
\item [] $\GNiceRes $: the pushforward of $\GNice$ to the residue
field of $\pi _{0}\Esource $; $\Esource$ is defined in equation
(\ref{bigtate}). 
\item [] $\Honda$: the Honda formal group law of height $n-1$ over $\k$.
\item [] $\Gabbrev$: A formal group law introduced in the proof of
Proposition~\ref{t-pr:E-phi-is-H} that is shown in that proof to be
the same as $\Honda$.
\end{enumerate}

We use the canonical orientation of $E (n)$, which provides a
coordinate so that 
\[
E (n)^{*} (\mathbf{C}P^{\infty }) = E
(n)^{*}\psb{x},
\]
 and as usual, if $\mu $ is the multiplication on
$\mathbf{C}P^{\infty }$, 
\begin{equation*}
\mu ^{*} (x) \in E (n)^{*}\psb{x,y} = E (n)^{*}\psb{x}\hat{\tensor{}} E
(n)^{*}\psb{y} = E (n)^{*} (\mathbf{C}P^{\infty }\times \mathbf{C}P^{\infty })
\end{equation*}
is a formal group law which we will denote $G (n) (x,y)$.

This formal group law has the feature that its $p$--series is given by 
$
\sum^{G (n)}_{i\leq n} v_{i}x^{p^{i}}.
$
Recall that
\begin{equation}\label{def:gltjw}
   \pi_{*}\gltjw{n} = \Zp \psb{u_{1},\dots ,u_{n-1}}[u,u^{-1}]
\end{equation}
with $|u_{i}|=0$ and $|u|=2$.  There is an isomorphism 
\[
\ELT{n}^{0} (\C P^{\infty }) \cong \pi_{0}\ELT{n}\psb{t}
\]
in terms of which the coproduct on $\ELT{n}^{0} (\C P^{\infty})$ is determined
by the formula 
\[
   t \mapsto \GLT{n} (s,t)
\]
where $\GLT{n}$ is the group law 
\[
   \GLT{n} (s,t) \eqdef u\GJW{n} (u^{-1}s,u^{-1}t)
\]
over $\pi_{0}\gltjw{n}$.

Since $v_{n-1}$ is a unit in 
$\pi_{*}\bigl(\Tate\EJW{n}\bigr)^{\wedge}_{I_{n-1}}$, we shall  
also consider the theory $\ERT$ obtained by adjoining an element $\w$
of degree $2$ such that  
\begin{equation} \label{eq:Root-of-vee-enn-minus-one}
    \w ^{q-1} = v_{n-1}
\end{equation}
and then completing with respect to the ideal $I_{n-1}= (p,\dots
,v_{n-2})$.  We prefer this choice to the usual parameter $u =
v_{n}^{1/ ( p^{n}-1)}$ (which gives $\gltjw{n}$) because the
functor $t$ will emphasize height $n-1$ behavior instead of height
$n$, and the normalization we choose to make things $2$--periodic leads
to simpler statements in section~\ref{sec:Isomorphism} and this section.
%

Equation \eqref{coefficient-calculation}  shows that there is an
isomorphism 
\[
  \pi_{*}\Tate \ERT \cong  \pi_{*}\ERT \lsb{\x}/\dps_{\GJW{n}} (\x).
\]
Proposition \ref{En-coeff-iso}  shows that $v_{n-1}$, and so also $\w$, is
a unit in the homotopy of the $I_{n-1}$--adic completion
\begin{equation}\label{bigtate}
    \Esource \eqdef \Bigl(\Tate \ERT\Bigr)\Iadic{n-1}
\end{equation}
which is thus $2$--periodic.  If $\GNice$ denotes the formal group law 
\begin{equation}\label{eq:GNice}
   \GNice (s,t) \eqdef \w \GJW{n} (\w^{-1}s,\w^{-1}t)
\end{equation}
(which implies $[p]_{F} (s) = w[p]_{G (n)} (w^{-1}s)$)
and we introduce elements 
\begin{align*}
     \w_{i} & = v_{i}\w^{- (p^{i}-1)} \\
     \y      & = \w \x
\end{align*}
of degree zero, then the argument of Proposition \ref{En-coeff-iso} 
shows that there are isomorphisms 
\begin{align}
   \pi_{0}\Esource  & \cong  \notag
\Bigl[\Zp[w_{1},\dots ,\w_{n-2},\w_{n}^{\pm 1}]\lsb{\y}/ 
               (\dps_{\GNice} (\y))\Bigr]\Iwadic \\
		    & \cong                       \label{eq:E-coefficients}
\Bigl[\Zp[\w_{1},\dots ,\w_{n-2}]\lsb{y}\Bigr]\Iwadic  \\
    \pi_{*}\Esource & \cong \Bigl[\Zp[\w_{1},\dots
,\w_{n-2}]\lsb{y}\Bigr]\Iwadic[\w,\w^{-1}]. \notag
\end{align}
The group law $\GNice$ is defined over 
$\pi_{0}E (n)[w]$  and hence over $\pi_{0}\Esource$; it is
$p$--typical, and its $p$--series satisfies the functional equation
\begin{equation}\label{eq:p-G-Nice}
   [p]_{\GNice } (t) = pt \fs{\GNice} \w_{1}t^{p} \fs{\GNice}
 \dots \fs{\GNice } \w_{n-2}t^{p^{n-2}} \fs{\GNice} t^{q} 
     \fs{\GNice} \w_{n} t^{pq}.
\end{equation}
We denote by $\GNiceRes $ the image of the group law $\GNice$ in the
residue field of $\pi _{0}TE$.

\begin{Proposition} \label{t-pr:p-series-residue-field}
The residue field of $\pi_{0}\Esource$ is $\k \lsb{\prootzero }$.   The element $\w_{n}$ of $\pi_{0}\Esource$ maps to $-\y^{(1-p)q}$
in the residue field.  The formal
group law $\GNiceRes$ has coefficients in the subring $\k[\prootzero^{-1}
]$, and its $p$--series satisfies the functional
equation  
\begin{equation}\label{eq:Function-eqn-G-Jack} 
  [p]_{\GNiceRes } (t) = t^{q}\fdiff{\GNiceRes } \prootzero^{(1-p)q} t^{pq}.
\end{equation}
\end{Proposition}

\begin{proof}
The statement about the residue field follows from equation
\eqref{eq:E-coefficients}.  Before we proceed further, recall that
since $\GNiceRes $ is $p$--typical and $p$ is odd, $[-1]_{\GNiceRes }
(t) = -t$.   Note that the image of equation
\eqref{eq:p-G-Nice}  in the residue field is 
\[
   [p]_{\GNiceRes } (t) = t^{q} \fs{\GNiceRes} \w_{n} t^{pq},
\]
so equation \eqref{eq:Function-eqn-G-Jack} follows from the assertion
that $\w_{n}$ maps to $-\y^{(1-p)q}$.  

Since $[p]_{\GNiceRes } (\y)=0$ we have
\begin{align*}
   0 &= \y^{q} \fs{\GNiceRes } \w_{n} \y^{pq} \\
\intertext{and so}
    \y^{q} & = [-1]_{\GNiceRes } (\w_{n} y^{pq}) = - \w_{n}\y^{pq}
\\
    \w_{n} &= -\y^{(1-p)q}.
\end{align*}

 Finally, $\GNiceRes $ is defined over the subring $\k[\y^{-1}]$
since $\GNice$ is actually defined over the polynomial ring in
$w_{1},\dots ,w_{n}$, and $(1-p)q$ is negative.
\end{proof}

So $\GNiceRes $ is a $p$--typical formal group law over $\k[\y^{-1}]$, of
height $n-1$ in the field $\k (\y)$ or $\k\lsb{\y}$.  On the other hand,
let $H$ be the Honda law of height $n-1$ over $\k$, characterized by
the fact that it is $p$--typical with $p$--series 
\[
  [p]_{H} (t) = t^{q}.
\]
Comparison with equation \eqref{eq:Function-eqn-G-Jack} shows that
\begin{equation}\label{eq:GNice-and-H-at-infty}
   \GNiceRes \equiv H \mod \y^{-1}.
\end{equation}

Now Lazard \cite{Lazard:LoisDeGroupes,Frohlich:FormalGroups} proves
that over the separable 
closure $\SepClos{\k (\y)}$, there is an
isomorphism of formal group laws 
\[
   \GNiceRes \cong H,
\]
which is in general not at
all canonical.  In this section we show that there is a
unique isomorphism $\Phi\co \GNiceRes \cong H$ which preserves
equation \eqref{eq:GNice-and-H-at-infty} in a suitable sense.

To make this precise, note that expanding a rational function as a
power series at infinity gives a map of fields
\[
     \k (\y) \rightarrow \k\lsb{\y^{-1}}
\]
which extends to 
\[
    \SepClos{\k (\y)} \rightarrow \SepClos{\k\lsb{\y^{-1}}}.
\]

\begin{Theorem} \label{t-th:Jacks-Psi}
There is a unique isomorphism $\Phi\co \E\to\Honda$ of formal group
laws over $\SepClos{\k (\y)}$ such that the image of $\Phi $ in $\SepClos{\k\lsb{\y^{-1}}}$
has coefficients in the image of 
\[
     \k\psb{\y^{-1}} \rightarrow \SepClos{\k\lsb{\y^{-1}}}.
\]
$\Phi $ gives an isomorphism over $\mathbf{F}_{p}\psb{y^{-1}}$ satisfying
\[
      \Phi (t) \equiv t \mod \y^{-1}.
\]
\end{Theorem}
We build up to the proof gradually; the proof itself appears after Proposition
\ref{t-pr:E-phi-is-H}.

\begin{Lemma} \label{t-le:tau}
There is a unique series $\tau (t) \in \k[\y^{-1}]\psb{t}$ such that
\[
    [p]_{\GNiceRes }  = [p]_{\Honda}\circ \tau.
\]
This power series has the properties that 
\begin{align*}
    \tau (t) & \equiv t \mod t^{2} \\
    \tau (t) & \equiv t \mod \y^{-1}. 
\end{align*}
\end{Lemma}

\begin{proof}
The equation which $\tau$ must satisfy is 
\[
[p]_{\GNiceRes } (t)  = \tau (t)^{q}.
\]
If 
\[
   \GNiceRes (s,t) = \sum_{i,j} b_{ij} s^{i}t^{j}
\]
then the functional equation \eqref{eq:Function-eqn-G-Jack} becomes
\begin{align*}
[p]_{\GNiceRes } (t) & = \sum_{i,j} b_{ij} 
                         (t^{q})^{i} (-\y^{(1-p)q}t^{pq})^{j} \\
                     & = \sum_{i,j} b_{ij} (-\y^{1-p})^{jq} t^{(i+pj)q}.
\end{align*}
So we must show that $b_{ij}$ has a unique $q$\th root.  If it does
then 
\[
  \tau (t) = \sum_{i,j} b_{ij}^{1/q} (-\y^{1-p})^{j}t^{(i+pj)}
\]
which shows that $\tau (t)\equiv t \mod t^{2}$.

If 
\[
   \GJW{n} (s,t) = \sum_{i,j} a_{ij} s^{i}t^{j}
\]
then by definition
\[
  \GNice (s,t) = \sum_{i,j} a_{ij}\w^{1-i-j} s^{i}t^{j},
\]
with $b_{ij} = a_{ij}\w^{1-i-j}$ homogeneous of degree zero.  The
$(p,v_{1},\dots ,v_{n-2})$ reduction
of $a_{ij}$ is of the 
form
\begin{equation} \label{eq:a-i-j}
   a_{ij} =  \sum_{a,b} c_{ab} v_{n-1}^{a}v_{n}^{b} \in \k[v_{n-1},v_{n}]
\end{equation}
for coefficients $c_{ab}$ depending on $i,j$.
Substituting 
\begin{align*}
     v_{n-1} & = \w^{q-1} \\
     v_{n}   & = \w^{pq-1} \w_{n} \\
             & = -\w^{pq-1} \y^{ (1-p)q}, 
\end{align*}
into equation \eqref{eq:a-i-j}, one has 
\[
 b_{ij}  = \w^{1-i-j} \sum_{a,b} 
           c_{ab} \w^{a (q-1)} (-\w^{pq-1}\y^{(1-p)q})^{b}.
\]
As $b_{ij}$ is homogeneous of degree zero, the exponent of $\w$ in
each term must add to zero, and so 
\begin{align*}
b_{ij}  & = \sum_{a,b} c_{ab} (-\y^{1-p})^{bq} \\
b_{ij}^{1/q} & = \sum_{a,b} c_{ab} (-\y^{1-p})^{b} 
\end{align*}
since $c_{ab} \in \mathbf{F}_{p}$ so $c_{ab}^{1/q} = c_{ab}$.
Thus one has 
\[
   \tau (t) = \sum_{i,j} r_{ij}t^{(i+pj)}, \text{ with}
\]
\[
  r_{ij} = \sum_{a,b} c_{ab} (-\y^{1-p})^{b+j}
\]
(recall the $c_{ab}$ depend on $i$ and $j$), which shows that 
$\tau (t) \equiv t \mod \y^{-1}$ as well.
\end{proof}

\begin{Proposition} \label{t-pr:Phi}
There is a unique power series 
$
\Phi (t) \in \k\psb{\y^{-1}}\psb{t}
$
such that 
\begin{align*}
 \Phi (t)  & \equiv t \mod t^{2};\\
 \Phi (t)  & \equiv t \mod \y^{-1}; \text{ and}\\
      [p]_{\Honda}\circ \Phi & = \Phi \circ [p]_{\GNiceRes }
\end{align*}
in $\k\psb{\y^{-1}}\psb{t}$.
\end{Proposition}

\begin{proof}
If $f$ in $A\psb{t}$ is a series with coefficients in a
ring $A$ of characteristic $p$, let $f^{\sigma}$ be the
corresponding series with coefficients those of $f$, raised to
the $q$\th power; thus 
$$
f^{\sigma}(t^{q})  =  (f(t))^{q} .
$$ 
In this notation, the equation supposedly satisfied by $\Phi$ takes 
the form 
\[
 \Phi^{\sigma}  =  \Phi \circ \tau^{\sigma}
\]
where 
$
\tau \in \k[\y^{-1}]\psb{t}
$
is the power series constructed in Lemma \ref{t-le:tau}.   
If $\tau^{-1}$ denotes the compositional inverse (not reciprocal) of
$\tau$, and similarly $\tau^{-\sigma}$ denotes $(\tau^{-1})^{\sigma}$, 
then this equation can be rewritten in the form 
\begin{align*}
\Phi & =  \Phi^{\sigma} \circ \tau^{-\sigma} \\
     & =  [ \Phi^{\sigma} \circ \tau^{-\sigma}]^{\sigma} 
          \circ \tau^{-\sigma} = \dots.
\end{align*}
If $\Phi$ is to be of the form 
\[
  \Phi (t)  \equiv t \mod \y^{-1},
\]
then  as $r$ grows, $\Phi^{\sigma^{r}}$ will converge to the
identity in the $(\y^{-1})$--adic topology, and we must have
\begin{equation}\label{eq:tau-to-phi}
  \Phi = 
   \lim_{r\to\infty} 
   \tau^{-\sigma^{r}}\circ \tau^{-\sigma^{r-1}}\circ 
	\dots \circ \tau^{-\sigma}.
\end{equation}
On the one hand, if the limit exists, then it certainly intertwines
the $p$--series.  On the other hand, the limit exists:
Lemma \ref{t-le:tau} asserts that 
\[
   \tau (t) \in  t + \y^{-1}\k[\y^{-1}]\psb{t},
\]
so $\tau^{-\sigma^{r}}$ converges to the identity in the
$(\y^{-1})$--adic topology.  The group of formal power series under
composition is complete with respect to the non-archimedean norm
defined by the degree of the leading term, so the infinite composite
\eqref{eq:tau-to-phi} converges because the sequence of composita 
converges to the identity.  It is easy to see in addition
that $\Phi$ inherits the property 
\[
 \Phi (t) \equiv t \mod t^{2}
\]
from $\tau$.
\end{proof}

\begin{Proposition} \label{t-pr:E-phi-is-H}
The power series $\Phi$ is the unique strict isomorphism 
\[
   \GNiceRes \xrightarrow{\Phi} \Honda
\]
of formal group laws over $\k\psb{\y^{-1}}$ such that 
\[
    \Phi (t) \equiv t \mod \y^{-1}.
\]
\end{Proposition}

\begin{proof}
Let us write $\Gabbrev $ for the formal group law
$\GNiceRes^{\Phi}$ where
\[
\GNiceRes^{\Phi} (x,y)=\Phi (\GNiceRes (\Phi ^{-1} (x),\Phi ^{-1} (y))).
\]
Since $[p]_{\Gabbrev } (t) = t^{q}$, the uniqueness of $\Phi
$ follows from Proposition \ref{t-pr:Phi}.

We need to show $\Gabbrev =\Honda.$
There is a canonical strict isomorphism
\[
   \Gtypical \xrightarrow{\rho} \Gabbrev,
\]
from a $p$--typical formal group law $\Gtypical$, defined over
$\k\psb{\y^{-1}}$.  Indeed $\rho$ is given by the formula 
\cite[A2.1.23]{Ravenel:GreenBook} 
\begin{equation}\label{eq:rho}
    \rho(t) = \sideset{}{^{\Gabbrev}}\sum_{p \notdivide r \geq 1}
                     [\mu (r)]_{\Gabbrev} [\tfrac{1}{r}]_{\Gabbrev}
		    \sideset{}{^{\Gabbrev}}\sum_{i=1}^{r} \zeta^{i}t,
\end{equation}
where 
\[
\mu (r) = \begin{cases} 0 & r \text{ is divisible by a square} \\
                       (-1)^{k}& r \text{ is the product of }k
                                   \text{ distinct primes}
\end{cases}
\]
and $\zeta $ is a primitive $r$th root of $1$.
We first claim that $\Gtypical=\Honda$; then we shall show that $\rho$
is the identity. 

As $\Gtypical$ and $\Honda$ are both are $p$--typical, it suffices to show that
$[p]_{\Gtypical} = [p]_{\Honda}$.  Since $\rho$ is a homomorphism of
groups, one has 
\begin{equation}\label{eq:Rho-conjugates-p}
     [p]_{\Gabbrev}\circ \rho = \rho\circ [p]_{\Gtypical}
\end{equation}
Using equation \eqref{eq:rho}, one has
\begin{align*}
[p]_{\Gabbrev}\circ \rho (t) & = 
                    \sideset{}{^{\Gabbrev}}\sum_{p \notdivide r \geq 1}
                     [\tfrac{\mu (r)}{r}]_{\Gabbrev}
		    \sideset{}{^{\Gabbrev}}\sum_{i=1}^{r} 
			[p]_{\Gabbrev} (\zeta^{i}t) \\
& =  \sideset{}{^{\Gabbrev}}\sum_{p \notdivide r \geq 1}
                     [\tfrac{\mu (r)}{r}]_{\Gabbrev}
		    \sideset{}{^{\Gabbrev}}\sum_{i=1}^{r} 
			\zeta^{iq}t^{q} \\
& = \rho(t^{q})
\end{align*}
since $q=p^{n-1}$ and $[p]_{\Gabbrev} (t) = t^{q}.$  Thus 
\[
   \rho([p]_{\Gtypical} (t)) = \rho(t^{q}) = \rho([p]_{\H} (t)).
\]

Next we must show that $\rho$ is the identity.  As
\[
    \Phi (t) \equiv t \mod \y^{-1}
\]
we have 
\[
     \GNiceRes \equiv \Gabbrev \mod \y^{-1}.
\]
But $\GNiceRes $ is $p$--typical, so
\[
     [\tfrac{1}{r}]_{\GNiceRes} 
     \sideset{}{^{\GNiceRes}}\sum_{i=1}^{r} \zeta^{i}t = 0
\]
for $p\notdivide r >1$; it follows from equation \eqref{eq:rho} that 
\begin{equation}\label{eq:rho-lead}
     \rho (t) \equiv t \mod \y^{-1}.
\end{equation}
But equation \eqref{eq:Rho-conjugates-p} and $[p]_{\Gabbrev } (t) = t^{q}$ imply that 
\[
     \rho (t^{q}) = \rho (t)^{q};
\]
if $\rho=\sum_{i\geq 1}\rho_{i}t^{i}$  then we must have 
\[
    \rho_{i} = \rho_{i}^{q}.
\]
Together with equation \eqref{eq:rho-lead} this implies $\rho_{1} =
1$.  For $i>1$ we have $\rho_{i} \in y^{-1}\k\psb{\y^{-1}}$, so
$\rho_{i}=0$. 
\end{proof}

\begin{proof}[Proof of Theorem \ref{t-th:Jacks-Psi}]
If 
$$\Psi(t) = \sum_{i \geq 0} \Psi_{i}t^{i+1}$$
is \emph{any} isomorphism from the group law $\GNiceRes$ to the group
law $\H$ over $\SepClos{\k (\y)}$, then it must satisfy the equation 
$$[p]_{\H} \circ \Psi  =  \Psi \circ [p]_{\E}.$$ 
By Lemma \ref{t-le:tau}, we must have 
\begin{equation}\label{eq:Psi-intertwines}
\Psi^{\sigma} = \Psi \circ \tau^{\sigma}.
\end{equation}
Because 
\[
   \tau (t) \equiv t \mod t^{2},
\]
the intertwining equation \eqref{eq:Psi-intertwines} can be rewritten
inductively as a sequence of  generalized Artin--Schreier equations 
\begin{equation}\label{define-psis}
\Psi_{i}^{q} - \Psi_{i}  =  \mbox{a polynomial in $\Psi_{j}$'s 
with $j<i$},
\end{equation}
beginning with 
\[
\Psi_{0}^{q}-\Psi_{0} = 0.
\]
Because of (\ref{define-psis}), the $\Psi_{i}$ are all
algebraic, and the Galois group 
of the extension they generate acts by translating the solutions by an element of the field with $q$ elements. 

The coefficients of $\Phi$ satisfy the same
equations in $\k\lsb{\y^{-1}}$.  Starting with $\Psi_{0}
= 1$, we may adjust each $\Psi_{i}$ by a Galois transformation so that
its image in $\SepClos{\k\lsb{\y^{-1}}}$ is $\Phi_{i}$.  The
resulting power series $\Psi$ is an isomorphism of formal group laws
in $\SepClos{\k (\y)}$, since it becomes one in
$\SepClos{\k\lsb{\y^{-1}}}$.

The uniqueness of $\Psi$ satisfying the hypotheses is a trivial
consequence of the uniqueness of $\Phi$ in Proposition \ref{t-pr:Phi}.
\end{proof}

It is the field $\k\lsb{\y}$ rather than $\k\lsb{\y^{-1}}$ which
appears in the Tate homology calculations. 
Expanding a rational function as a power series at zero gives an
embedding $\k (y)\to\k\lsb{y}$, which extends to an embedding
\[
   \SepClos{\k (y)} \rightarrow \SepClos{\k\lsb{y}}.
\]
Thus we have 
\begin{Corollary} \label{t-co:Phi-at-last}
There is a unique strict isomorphism  $\Phi\co \E \rightarrow\H$,
satisfying
\begin{itemize}
\item[\rm(1)] the coefficients of $\Phi $ are in the subfield 
$\SepClos{\k (\y)}$ and 
\item[\rm(2)] the expansion of $\Phi $ at $y=\infty $ is a power 
series with coefficients in $\k\psb{\y^{-1}}$, congruent to the
identity modulo $\y^{-1}$. 
\end{itemize} 
\end{Corollary}

\section{A map of ring spectra}
\label{sec:Generalized-chern-character}

\subsection{Lubin and Tate's theory of lifts} \label{sec:LubinTate}

We recall briefly the deformation theory of Lubin and Tate
\cite{LubinTate:FormalModuli}.  Suppose that $\F$ is a field of characteristic
$p$.    

\begin{Definition} \label{def:LiftsOfField}
A \emph{lift} of $\F$ is a pair $(A,i)$
consisting of 
\begin{itemize}
\item[(1)] a Noetherian complete local ring $A$ with residue field $A_{0}$;
\item[(2)] a map of fields $i\co \F\to A_{0}$.
\end{itemize}
A map $f\co (A,i)\to (B,j)$ (or $(A,i)$--\emph{algebra}) is a local
homomorphism  
\[
     A \xrightarrow{f} B
\]
such that $j=f_{0}\circ i$, where $f_{0}\co A_{0}\rightarrow B_{0}$ is
induced by $f$. 
\end{Definition}    
We shall abbreviate $(A,i)$ to $A$ when $i$ is clear from
context.  

Suppose that $\Gamma$ is a formal group law of finite height $n$ over a 
field $\F$ of characteristic $p$, that $(A,i)$ is a lift of $\F$,
and that $(B,j)$ is an $(A,i)$--algebra.

\begin{Definition} \label{def:deformation}
A \emph{deformation of } $\Gamma$ \emph{ to }$(B,j)$
is a pair $(G,\phi)$ 
consisting of 
\begin{itemize}
\item[(1)] a formal group law $G$ over $B$;
\item[(2)] an isomorphism of group laws $\phi\co j^{*}\Gamma \approx
G_{0}$, where $G_{0}$ denotes the group law over $B_{0}$ induced by
$G$. 
\end{itemize}
Two deformations $(G,\phi)$ and $(G',\phi')$ are \st--isomorphic if
there is an isomorphism $c\co G\to G'$ such that 
\[
       \phi'=c_{0}\circ \phi\co j^{*}\Gamma \rightarrow G'_{0}
\]
in $B_{0}$.
\end{Definition}

The set of \st--isomorphism classes of deformations to $(B,j)$ is a
functor from $(A,i)$--algebras 
to sets;  
Lubin and Tate show that this functor is representable.  Namely, 
let
\begin{align*}
   R  & = A\psb{u_{1},\dots ,u_{n-1}},\\
  C_{p^{i}} (x,y) & = \frac{1}{p}[(x+y)^{p^{i}}-x^{p^{i}}-y^{p^{i}}],
\end{align*}
and let $(G,\phi)$ be any deformation of $\Gamma$ to $R$ such
that 
\begin{equation} \label{eq:Lubin-Tate-criterion}
\begin{split} 
G_{0} & = \Gamma \\
G (s,t)&  \equiv s + t + u_{i}C_{p^{i}} (x,y) \mod u_{1},\dots
,u_{i-1},\text{ and degree }p^{i}+1
\end{split}
\end{equation}
for $\epsilon$ a unit of $A_{0}$ and $1\leq i\leq n-1$.
Lubin and Tate show that such deformations exist; we shall call such a
group law a \emph{Lubin-Tate lift} of $\Gamma$.
Theorem 3.1 of \cite{LubinTate:FormalModuli} may be phrased as
follows.

\begin{Theorem}\label{t-th:Lubin-Tate}
If $(G',\phi')$ is a deformation of $\Gamma$ to an $(A,i)$--algebra $(B,j)$,
then there is a unique map of $(A,i)$--algebras $f\co R\to (B,j)$ such
that there is a \st--isomorphism 
\[
     (f^{*}G,f_{0}^{*}\phi) \xrightarrow{c} (G',\phi').
\]
Moreover the \st--isomorphism is unique.
\end{Theorem}

\subsection{The map of ring theories}

A classical theorem of I Cohen (see
\cite{Schoeller:GroupesCorpsNonParfait}) asserts that for 
any field $\F$ of characteristic $p$, there is an essentially unique
complete discrete valuation ring $C_{\F}$ with maximal ideal generated
by $p$ and residue field $\F$; in particular, $C_{\F}$ is Noetherian.
If $\F$ is perfect, then $C_{\F}$ is just 
the ring $\mathbf{W}_{\F}$ of Witt vectors of $\F$, but in general it is a
subring of the Witt vectors, and (although it is a functor) it is
not very easily described. However, when the degree $[\F:\F^{p}]$
is finite (instead of $1$ as in the perfect case), eg if $\F =
\k\lsb{\y}$, then the Cohen ring is still relatively tractable.
In this case, for example, we have

\begin{Lemma} \label{t-le:CohenOfRes}
The Cohen ring of $\k\lsb{\y}$ is isomorphic to the $p$--adic
completion of $\Zp\lsb{\y}$.  Any such isomorphism extends
to an isomorphism 
\[
  \C_{\k\lsb{\y}}\psb{\w_{1},\dots \w_{n-2}} \cong \pi_{0}\Esource.
\]
\end{Lemma}
\begin{proof}
The $p$--adic completion of $\Zp\lsb{\y}$ is a complete discrete
valuation ring with maximal ideal generated by $p$ and residue field
$\k\lsb{\y}$, so it is isomorphic to the Cohen ring of
$\k\lsb{\y}$.  Indeed Schoeller \cite[\S
8]{Schoeller:GroupesCorpsNonParfait} shows that an isomorphism is given
by a choice of $p$--base for $\k\lsb{\y}$ and representatives of that
$p$--base in $(\Zp\lsb{\y})^{\wedge}_{p}$.

Equation \eqref{eq:E-coefficients} reduces the second part of the lemma
to the observation that the rings 
\[
  \bigl(\Zp\lsb{\y}\bigr)^{\wedge}_{p}\psb{\w_{1},\dots ,\w_{n-2}} \cong
  \bigl(\Zp\psb{\w_{1},\dots,\w_{n-2}} \lsb{\y} \bigr)\Iwadic
\]
are isomorphic.
\end{proof}

\begin{Proposition} \label{t-pr:LubinTate-for-Tate}
The formal group law $\GNice$ over $\pi _{0}\Esource$ is a Lubin--Tate
lift of the group law $\GNiceRes $ of height $n-1$ over $\k
\lsb{\prootzero }$ to $\C_{\k\lsb{\y}}$--algebras. 
\end{Proposition}

\begin{proof}
It is a standard fact about $\GJW{n}$ that 
\[
G (s,t)  \equiv s + t + v_{i}C_{p^{i}} (s,t) \mod v_{1},\dots
,v_{i-1},\text{ and degree }p^{i}+1
\]
for $1\leq i\leq n-1$.  By the definition \eqref{eq:GNice} of
$\GNice$, one then has
\begin{align*}
\GNice (s,t)&  \equiv s + t + w_{i}C_{p^{i}} (s,t) \mod w_{1},\dots
,w_{i-1},\text{ and degree }p^{i}+1
\end{align*}
for $1 \leq i \leq n-2$.  It
follows that $\GNice$ satisfies equations 
\eqref{eq:Lubin-Tate-criterion} for the group law $\Gamma =\GNiceRes$.
\end{proof}

The formal group law $\GLT{n-1}$ over $\pi _{0}\gltjw{n-1}$ has image
$\Honda $ in  $\k$, so by Corollary \ref{t-co:Phi-at-last} the pair
$(\C_{\khuge}\hat{\otimes}\GLT{n-1}, \Phi)$ is a deformation of
$\GNiceRes $ to the $\C_{\k\lsb{\y}}$--algebra 
\[
\C_{\khuge}\hat{\tensor{\Zp }}\pi_{0}\gltjw{n-1} \cong
\C_{\khuge}\psb{u_{1},\dots ,u_{n-2}}.
\]
Theorem
\ref{t-th:Lubin-Tate} provides a ring homomorphism 
\begin{equation}\label{eq:tE-to-Enmo}
     \pi_{0}\Esource \xrightarrow{f}
\C_{\khuge}\hat{\tensor{\Zp }}\pi_{0}\gltjw{n-1}
\end{equation}
and an isomorphism of formal group laws
\begin{equation}\label{eq:GNice-to-GLT}
     f^{*}\GNice \xrightarrow[\cong]{c} \C_{\khuge}\hat{\tensor{\Zp }}
\GLT{n-1}
\end{equation}
such that $c_{0} = \Phi$.

\begin{Proposition} \label{t-pr:tE-to-CEnmo-iso-after-extension}
The map 
\[
\C_{\khuge}\hat{\otimes} f\co \C_{\khuge}\hat{\otimes}_{\C_{\k \lsb{y}}}\pi_{0}\Esource
\rightarrow \C_{\khuge}\hat{\tensor{\Zp }}\pi_{0}\gltjw{n-1}
\]
is an isomorphism.
\end{Proposition}

\begin{proof}
The ring on the left represents deformations of
$\khuge\hat{\otimes}\GNiceRes$ to\break
$\C_{\khuge }$--algebras, while the ring on
the right represents deformations of\break
$\khuge \hat{\otimes} H$ to
$\C_{\khuge}$--algebras.  The isomorphism $\Phi$ of Corollary
\ref{t-co:Phi-at-last} induces an isomorphism between these functors. 
\end{proof}

There are isomorphisms 
\begin{align*}
     \Esource ^{0} (\C P^{\infty} ) &\cong \Esource ^{0}\psb{t_{1}} \\
     \ELT{n-1}^{0} (\C P^{\infty} ) &\cong \ELT{n-1}^{0}\psb{t_{2}} \\
\end{align*}
such that the  formal group law $\GNice$ expresses the coproduct on $\Esource
^{0} (\C P^{\infty} )$, and the formal group law $\GLT{n-1}$ expresses the
coproduct on $\ELT{n-1}$.   
A standard argument \cite{Miller:Ell} using Landweber's exact functor
theorem gives:

\begin{Theorem} \label{t-th:Ring-map}
There is a canonical map of ring theories 
\[
\Esource \xrightarrow{\Theta} \C_{\khuge}\hat{\tensor{\Zp }}\gltjw{n-1}
\]
whose value on coefficients is $f$, and whose value on 
$\Esource^{0} (\C P^{\infty} )$ is determined by the equation 
\[
  \Theta (t_{1}) = c(t_{2}).
\]
Here $c$ is the isomorphism \eqref{eq:GNice-to-GLT}.
\qed
\end{Theorem}

\begin{Corollary}\label{iso-ring-spec}
$\C_{\khuge}\hat{\tensor{\Zp }} \Theta $ gives an equivalence of
complex oriented ring spectra
\[
\C_{\khuge}\hat{\tensor{\C_{\k\lsb{y}}}}TE \rightarrow
\C_{\khuge}\hat{\tensor{\Zp }} \gltjw{n-1}.
\]
\end{Corollary}

Composing the map from Theorem~\ref{t-th:Ring-map} with the maps 
\[
  \EJW{n} \rightarrow \ERT \rightarrow \Esource 
\]
gives a canonical map of ring theories 
\[
 \EJW{n} \rightarrow  \C_{\khuge}\hat{\tensor{\Zp }}\gltjw{n-1},
\]
our generalized Chern character.

\end{document}